\documentclass[12pt]{article}

\usepackage{graphicx}

\usepackage[reqno,tbtags]{amsmath}

\usepackage{amssymb}

\def\thebibliography#1{\vspace{0.5cm} {\flushleft \bf
References \\ } \list
 {\arabic{enumi}.}{\settowidth\labelwidth{[#1]}\leftmargin\labelwidth
 \advance\leftmargin\labelsep
 \usecounter{enumi}}
 \def\newblock{\hskip .11em plus .33em minus -.07em}
 \sloppy
 \sfcode`\.=1000\relax}

\protect

\protect

\protect

\protect

\newskip\halflineskip

\newtheorem{theorem}{Theorem}[section]
\newtheorem{corollary}[theorem]{Corollary}
\newtheorem{proposition}[theorem]{Proposition}
\newtheorem{lemma}[theorem]{Lemma}
\newtheorem{definition}[theorem]{Definition}

\newtheorem{remark}[theorem]{Remark}

\newenvironment{proof}{\noindent{\bf Proof. }\rm}
{\unskip\nobreak\hfil\penalty50\hskip1em\hbox{}
\nobreak\hfill\qed\par\smallskip}
\def\qed{\vrule height1ex width1ex depth0pt}

{\protect\nopagebreak\section*{\protect\raggedright \small \it
Acknowledgements}}%
{}

{\protect\nopagebreak\section*{\protect\raggedright \normalsize
Acknowledgement}}%
{}

\protect

\protect\newcommand{\sectionhead}[1]%
{\section {#1}}

\def\abstract{\vspace{0cm} \noindent
{\bf \footnotesize Abstract. }\footnotesize }

\title{{\bf  Monotone operator functions, gaps and power moment problem}}

\author{{\small {\sc HIROYUKI OSAKA }} \\
{\small {\it Department of Mathematical Sciences, Ritsumeikan University,}}\\
{\small {\it Kusatsu, Shiga 525-8577, Japan} } \\
{\small {\it e-mail: osaka@se.ritsumei.ac.jp} } {\small {\it FAX:
+81 77 561 2657 \ \ tel: +81 77 561 2656 \vspace{0.5cm}}}
\\
{\small {\sc SERGEI SILVESTROV}}
\\
{\small {\it Centre for Mathematical Sciences,} }\\
{\small {\it Lund University, Box 118, SE-22100 Lund, Sweden.} }\\
{\small {\it e-mail: sergei.silvestrov@math.lth.se } }
{\small{\it FAX:  +46 46 2224010 \ \ tel: +46 46 2228854 \vspace{0.5cm} }}\\
{\small {\sc JUN TOMIYAMA }} \\
{\small{\it Prof. Emeritus of Tokyo Metropolitan University,}}\\
{\small{\it 201 11-10 Nakane 1-chome, Meguro-ku, Tokyo, Japan }}\\
{\small{\it e-mail: jtomiyama@fc.jwu.ac.jp }} }
\date{March 1, 2006}

\begin{document}
\maketitle

\footnotetext{Dedicated to the memory of late Professor Gert K.
Pedersen.} \footnotetext{ This work was supported by The Swedish
Foundation for International Cooperation in Research and Higher
Education (STINT), Crafoord Foundation, The Royal Physiographic
Society in Lund and Open Research Center Project for Private
Japanese Universities: matching fund from MEXT, 2004-2008.}


 \pagebreak

\begin{abstract}
The article is devoted to investigation of the classes of
functions belonging to the gaps between classes $P_{n+1}(I)$ and $P_{n}(I)$ of matrix monotone functions for full matrix algebras of successive dimensions. In this paper we address the problem of characterizing polynomials belonging to the gaps $P_{n}(I) \setminus P_{n+1}(I)$ for bounded intervals $I$. We show that solution of this problem is closely linked to solution of truncated moment problems, Hankel matrices and Hankel extensions.
Namely, we show that using the solutions to truncated moment
problems we can construct continuum many polynomials in the gaps.
We also provide via several examples some first insights into the
further problem of description of polynomials in the gaps that are
not coming from the truncated moment problem. Also, in this
article, we deepen further in another way into the structure of
the classes of matrix monotone functions and of the gaps between
them by considering the problem of position in the gaps of certain
interesting subclasses of matrix monotone functions that appeared
in connection to interpolation of spaces and in a prove of the
L{\"o}wner theorem on integral representation of operator monotone
functions.

\noindent \thanks{\footnotesize {\bf  Keywords}:  operator monotone functions, matrix monotone functions} \\
\noindent \thanks{ \footnotesize {\bf  Mathematics Subject
Classification 2000:} Primary 26A48 Secondary  47A56, 47A63,
44A60}
\end{abstract}

\section{Introduction.} \label{sec:intr}
A real-valued continuous function $f:I\rightarrow \mathbb{R}$ is said to be
matrix monotone of order $n$ over an interval
$I$, if \begin{equation} \label{ineq:monot}
x\leq y \quad \Rightarrow \quad f(x)  \leq f(y)
\end{equation}
for any two self-adjoint $n \times n$ matrices $x$ and $y$ with eigenvalues in $I$. We denote the class of all such functions by $P_n(I)$.
A real-valued continuous function $f: I \mapsto \mathbb{R}$ on a
(non trivial) interval $I \neq \mathbb{R}$ is called
operator monotone if
the implication \eqref{ineq:monot} holds for any pair of bounded operators
$x, y \in B(H)$ on an infinite-dimensional separable
Hilbert space $H$ with their spectra in $I$.
We denote the class of all operator monotone functions
over an interval $I$ by $P_{\infty}(I)$, or simply by $P_{\infty}$ when the choice of the interval is clear from context.
For each positive integer $n$,
the proper inclusion $P_{n+1}(I) \subsetneq P_{n}(I)$ holds.
This fact has been stated in \cite{Donoghuebook}, but the complete proof of this appeared first in \cite{HansenJiTomiyama-art}.
The gaps $P_{n}(I)\setminus P_{n+1}(I)$ between classes of monotone matrix functions were also recently addressed in
\cite{OsakaSilvestrovTomiyamaIJM-art} and \cite{Nayak1}.
For infinite-dimensional Hilbert space, the set of operator monotone functions on $I$ can be shown to coincide with the intersection $$P_{\infty}(I) = \bigcap_{n=1}^{\infty} P_n(I),$$ or in other words a function is operator monotone if and only if it is matrix monotone of order $n$ for all positive integers $n$ \cite[Chap.5, Proposition 5.1.5 (1)]{HiaiYanagibook}.

The proof of non-emptiness of gaps 
$P_{n}(I)\setminus P_{n+1}(I)$
in \cite{HansenJiTomiyama-art} is constructive, by exhibiting for each positive integer $n$ an explicit function in the gap.
Moreover, for any bounded interval and each positive integer $n$, that function in the gap $P_{n}(I)\setminus P_{n+1}(I)$ exhibited in \cite{HansenJiTomiyama-art} was a polynomial, thus suggesting that there might be more polynomials in the gaps for any bounded interval, thus leading directly to an interesting problem of characterizing such polynomials. For the unbounded interval $(0,+\infty)$ it can be shown that there are no polynomials in the gaps. However, the unbounded interval can be bijectively mapped onto a bounded interval using an operator monotone fractional M{\"o}bius transformation with operator monotone inverse, and then any polynomial in the gap over that bounded interval, after proper composition with those fractional M{\"o}bius transformations, yields a rational function from the gap over the unbounded interval.

In this paper we address the problem of characterizing polynomials belonging to the gaps $P_{n+1}(I) \subsetneq P_{n}(I)$ for bounded intervals $I$. 
We show that solution of this problem is closely
linked to solution of truncated moment problems, Hankel matrices
and Hankel extensions. Namely, we show that using the solutions to truncated moment problems we can construct continuum many
polynomials in the gaps. We also provide via several examples some first insights into the further problem of description of
polynomials in the gaps that are not coming from the truncated
moment problem.

Also, in this article, we deepen further in another way into the
structure of the classes $P_{n}$ and of the gaps by considering a certain interesting subclass of functions inside $P_{n}$. 
This class of functions, denoted by ${\cal M}_n((0,+\infty))$, has been defined in \cite{SparrMathScandart}, 
as consisting of real-valued functions $h$ on $(0,\infty)$ such that for 
$a_j \in \mathbb{R}$, $\lambda_j > 0$ and $j = 1,\dots, 2n$ the following implication holds:
\begin{equation} \label{rel:defMnfunclas}
\left(\sum_{j=1}^{2n} a_j \frac{t \lambda_j - 1}{t+\lambda_j} \geq 0 \mbox{ for } t > 0, \sum_{j=1}^{2n} a_j = 0 \right) \Rightarrow
\left(\quad \sum_{j=1}^{2n} a_j h(\lambda_j) \geq 0\right).
\end{equation}
It was shown in \cite{SparrMathScandart} that $$P_{n+1}((0,+\infty)) \subseteq {\cal M}_n((0,+\infty)) \subseteq P_{n}((0,+\infty))$$
for any positive integer $n$, and so $
P_{\infty}=\bigcap_{n=1}^{\infty} P_{n}((0,+\infty)) =
\bigcap_{n=1}^{\infty} {\cal M}_{n}((0,+\infty)).$ In
\cite{SparrMathScandart}, an explicit example, showing that $P_2
\setminus {\cal M}_2 \neq \emptyset$, has been pointed out, thus
particularly implying that $P_2((0,+\infty)) \setminus
P_3((0,+\infty)) \neq \emptyset$. Proving that $P_n((0,+\infty))
\setminus {\cal M}_n((0,+\infty)) \neq \emptyset$ and ${\cal
M}_n((0,+\infty)) \setminus P_{n+1}((0,+\infty)) \neq \emptyset$
for an arbitrary $n$ is still an open problem. The unbounded
interval $(0,+\infty)$ is a union of inclusion increasing set of
bounded intervals $(0,+\infty) = \cup_{a>0} (0,a)$. In this
article we consider the classes of functions ${\cal M}_n(I)$ on
the bounded intervals. The definition is the same up to just
replacing $(0,+\infty)$ by the bounded interval $I$. The content
of the class ${\cal M}_n(I)$ differs from ${\cal
M}_n((0,+\infty))$. However, we provide in this article a proof that the  inclusions $$P_{n+1}(I) \subseteq {\cal M}_n(I) \subseteq P_{n}(I)$$ hold even for any bounded interval of the form $(0,a)$ or $(0,a]$ and all positive integers $n$. Therefore, we can conclude that $\bigcap_{n=1}^{\infty} {\cal M}_{n}(I)= P_{\infty}(I)$. The problem of proving or disproving the existence of the non-empty gap $P_{n}(I) \setminus {\cal M}_n(I)$ is also an open problem both for the bounded interval $I$ and for $(0,+\infty)$. However, while the example of function in the gap $P_{2}((0,+\infty)) \setminus {\cal M}_{2}((0,+\infty))$
constructed in \cite{SparrMathScandart} is non-polynomial due to
lack of polynomials and also seems to be difficult to extend to an arbitrary $n$, in the case of the bounded interval, we show in this article how to construct explicitly infinitely many
polynomials in the gap $P_{n}(I) \setminus P_{n+1}(I)$ for any
$n$. Thus a natural problem is to describe position of these
polynomials with respect to the gaps $P_{n}(I) \setminus {\cal
M}_n(I)$ and ${\cal M}_n(I) \setminus P_{n+1}(I) $. We succeeded
to investigate this problem for the polynomial in the gap
constructed in \cite{HansenJiTomiyama-art} for $n=2,3,4,5$, using numerical computations in Maple.

\section{Polynomial Monotone Matrix Functions.}
\label{sec:PolyGaps}

\begin{proposition}
\label{th:nopolyPinfint} The only polynomials belonging to the
class $P_l([0,+\infty))$ for an integer $l>1$ are polynomials of
the form $a t + b$ where $a \geq 0$.
\end{proposition}

\begin{proof}
Let $p_n(t) = \sum_{j=0}^na_jt^{n-j}$ be a polynomial in $P_l([0, +\infty))$
with $a_0 \not= 0$.
Then for any $C, D \in M_l$ such that $0 \leq C \leq D$ and any $\lambda > 0$ we have $0 \leq \lambda C \leq \lambda D$ 
and hence
\begin{align*}
a_n I &\leq \sum_{j=0}^na_j(\lambda C)^{n-j}
\leq \sum_{j=0}^n(\lambda D)^{n-j}\\
\frac{a_n}{\lambda^n}I &\leq \sum_{j=0}^n\lambda^{-j}C^{n-j}
\leq \sum_{j=0}^na_j\lambda^{-j}D^{n-j}
\end{align*}
which, after passing to the limit $\lambda \rightarrow + \infty$,
yelds
$0 \leq a_0C^n \leq a_0D^n$
implying $0 < a_0$ and $0 \leq C^n \leq D^n$.
This holds for arbitrary choice of $0 \leq C \leq D$ only if $n = 1$,
since $f(t) = t^n \notin P_l([0,+\infty)) \subset P_2([0,+\infty))$ when
$n > 1$ and $l \geq 2$. Thus, $p_n(t) = a_{n-1}t + a_n$.
When $a_{n-1} \not= 0$, from the same argument as in the first we have
$a_{n-1} > 0$, which is exactly what had to be proved.
\end{proof}

The situation is totally different on the finite intervals. There
polynomials of high degree than one can be matrix monotone of
order $n$. There is no contradiction here since the
transformations between a finite and an infinite interval do not
map polynomials into polynomials. Usually a M{\"o}bius
transformation can be used for this purpose, and in this case the
polynomial on a finite interval will be transformed into a
rational function on an infinite interval.

Let $g_n (t)= t + \frac{1}{ 3} t^3 + \dots + \frac{1}{2n-1}
t^{2n-1}$ , where $n$ is some positive integer. In
\cite{HansenJiTomiyama-art} it was proved that there exists
$\alpha_n > 0$ such that $g_n \in P_n([0,\alpha_n)) \setminus
P_{n+1}([0,\alpha_n))$, and consequently $f_n  = g_n \circ h_{n}
\in P_n \setminus P_{n+1}$, where $h_n(t)$ is the M{\"o}bius
transformation $h_n(t)= \frac{\alpha_n t}{1+t}$, operator monotone
on $[0,\infty)$, with the inverse $h_n^{\circ (-1)}(t)=
\frac{t}{\alpha_n-t}$ operator monotone on $[0,\alpha_n)$.

Note that two compact intervals can be however mapped to each
other with some polynomial of degree one
$\alpha t + \beta$ with $\alpha > 0$, an operator monotone function on any interval.
Namely, the bounded
interval with end points $u_1<v_1$ is
mapped to the bounded interval with end points $u_2<v_2$
by the map
$h(t)=\frac{v_2-u_2}{v_1-u_1}t +
\frac{u_2v_1-v_2 u_1}{v_1-u_1}$,
with the composition inverse
$$h^{\circ (-1)}(t) = \frac{v_1-u_1}{v_2-u_2}t-\frac{u_2v_1-v_2 u_1}{v_2-u_2}$$
which are both operator monotone since $\frac{v_2-u_2}{v_1-u_1}>0$
and $\frac{v_1-u_1}{v_2-u_2}>0$. The type of the interval with
respect to the inclusion or exclusion of the end points is
preserved by this map. Moreover, this map transforms polynomials
matrix monotone of order $n$ on one interval into polynomials of
the same degree and matrix monotone of order $n$ on the other
interval. In particular, $[0,a)$ is transformed to $[u,v)$ by the
map $h(t)=\frac{v-u}{a}t +u$ with the composition inverse
$h^{\circ (-1)}(t) = \frac{a}{v-u}t-\frac{au}{v-u}$. The interval
$[0,a)$ is mapped to the interval $[0,b)$ by the map
$h(t)=\frac{v}{a}t +u$ with the composition inverse $h^{\circ
(-1)}(t) = \frac{a}{v}t-\frac{au}{v}$. The interval $[-1,1]$ is
mapped to the interval $[u,v]$ by the map $h(t)=\frac{v-u}{2}t
+\frac{v+u}{2}$ with the composition inverse $h^{\circ (-1)}(t) =
\frac{2}{v-u}t-\frac{v+u}{v-u}$. The interval $[-1,1]$ is mapped
to the interval $[0,a]$ by the map $h(t)=\frac{a}{2}t
+\frac{a}{2}$ with the composition inverse $h^{\circ (-1)}(t) =
\frac{2}{a}t-1$. Keeping these considerations on maps of the
intervals in mind, we will work on the intervals containing $0$ or
other intervals convenient for the proofs, making clear from our
statements or by specially pointing out when the choice of the
interval is not essential.

We will make use of the following conditions concerned with
$n$-monotonicity of functions on an interval \cite{Donoghuebook},
restricting formulation to the functions which are infinitely
differentiable, which is suited to our considerations. For every
such function and every positive integer $n$ define the matrix
$M_n(f;t)
=\left(\frac{f^{(i+j-1)}(t)}{(i+j-1)!}\right)_{i,j=1}^n$. If $f\in
P_n((a,b))$ for $n\geq 2$, then $M_n(f;t)\geq 0$ and
$f^{(2n-3)}(t)$ is convex on $(a,b)$, by \cite[Theorem VI, Ch.
VII]{Donoghuebook}. Conversely, if $M_n(f;t)\geq 0$ and the
derivative $f^{(2n-3)}(t)$ is positive and convex, then $f\in
P_n((a,b))$, by \cite[Theorem V, Ch. VIII]{Donoghuebook}).

\begin{theorem} \label{the:polynomials}
Let $I\subset \mathbb{R}$ be a bounded interval on the real line.
There are no polynomials of degree $1 < \deg(f) < 2n-1$ in the
class $P_n(I)$, and there exists a polynomial $f$ of any order
$\deg(f) \geq 2n-1$ in $P_n(I)$. Any polynomial of degree $\deg
(f) = 2n-1$ or $\deg (f) = 2n$ belonging to $ P_n(I)$ lies in the
gap $f\in P_n(I) \setminus P_{n+1}(I)$.
\end{theorem}

\begin{proof}
Let $f(t)= c + \sum_{j=0}^{k-1} b_{j} t^{j+1}$, where $b_{k-1}\neq
0$ and $1< k=\deg(f) < 2n-1$. We consider two cases, of odd and
even $k$. Let $k=2l$ with $l\geq 1$. Then
$$
M_{l+1}(f;0) = \left(\begin{array}[c]{lllll}
b_0 & b_1 & \dots & b_{l-1} & b_{l} \\
b_1 & \dots & \dots & b_{l} & b_{l+1} \\
\vdots & \vdots &\vdots & \vdots & \vdots \\
b_{l-1} & b_{l}& \dots & b_{2l-2} & b_{2l-1} \\
b_{l}    & b_{l+1} & \dots & b_{2l-1} & 0
\end{array}
\right)$$ Since $k< 2n-1$, the matrix $M_{l+1}(f;0)$ is contained
as the principle upper left conner submatrix of $M_n(f;0)$.
Because
$$\det \left(\begin{array}[c]{ll}
b_{k-2} & b_{k-1} \\
b_{k-1} & 0
\end{array}
\right)
= -(b_{k-1})^2 <0,
$$
the matrix $M_n(f;0)$ is not positive definite and therefore $f
\not\in P_n(I)$. In the odd case, that is for $k=2l-1$, $l\geq 2$
and $b_{2l-2}\neq 0$, one has
$$
M_{l+1}(f;0) = \left(\begin{array}[c]{llllll}
b_0 & b_1 & \dots & b_{l-2} & b_{l-1} & b_{l} \\
b_1 & b_2 & \dots & b_{l-1} & b_{l} & b_{l+1} \\
\vdots & \vdots &\vdots & \vdots & \vdots & \vdots\\
b_{l-2} & b_{l-1}& \dots & b_{2l-4} &b_{2l-3} & b_{2l-2} \\
b_{l-1} & b_{l}& \dots & b_{2l-3} & b_{2l-2} & 0 \\
b_{l}    & b_{l+1} & \dots & b_{2l-2} & 0        & 0
\end{array} \right).
$$
Suppose that $f \in P_n(I)$ and hence $M_n(f;0) \geq 0$.
Then $b_{k-1}=b_{2l-2} >0$, since $M_n(f;0) \geq 0$ and
since $b_{k-1}=b_{2l-2}\neq 0$ as the highest coefficient
of the polynomial $f$.
Thus
$$
\det \left(\begin{array}[c]{lll}
b_{2l-4} & b_{2l-3} & b_{2l-2} \\
 b_{2l-3} & b_{2l-2} & 0 \\
 b_{2l-2} & 0 & 0 \end{array}
\right)
= -(b_{2l-2})^3 <0.
$$
Hence the matrix $M_n(f;0)$ is not positive semi-definite which
contradicts to the assumption $f \in P_n(I)$. Therefore $f \not\in
P_n(I)$.

In \cite{HansenJiTomiyama-art} it was proved that for any positive
integer $n$ there exists $\alpha_n > 0$ such that $g_n (t)= t +
\frac{1}{ 3} t^3 + \dots + \frac{1}{2n-1} t^{2n-1} \in
P_n([0,\alpha_n)) \setminus P_{n+1}([0,\alpha_n))$. Consequently,
if $I=[u,v)$, then  $\tilde{g}_n  = g_n \circ h^{\circ (-1)} \in
P_n([u,v)) \setminus P_{n+1}([u,v))$ where $\tilde{g}_n$ is the
polynomial of degree $2n-1$ obtained by composition of $g_n$ with
the operator monotone affine transformation $h^{\circ (-1)}(t) =
\frac{\alpha_n}{v-u}t-\frac{\alpha_n u}{v-u}$, mapping interval
$[u,v)$ onto $[0,\alpha_n)$. In order to show existence of the
polynomials of the even degree in the gap, take $p_n (t)= t +
\frac{1}{ 3} t^3 + \dots + \frac{1}{2n-1} t^{2n-1} + a t^{2n}$. By
the first statement of the theorem we have already proved that
$p_n \not\in P_{n+1}([0,\alpha))$ for any $\alpha >0$ since $\deg
(p_n) = 2n <2(n+1)-1=2n+1$. Since ${\rm det}M_n(p_n;0) = {\rm
det}M_n(g_n;0)$, there exists $\alpha'_n > 0$ such that $p_n \in
P_n([0,\alpha'_n))\backslash P_{n+1}([0,\alpha'_n))$. Therefore,
the polynomial $p_n \circ h^{\circ(-1)}$ of degree $2n$ belongs to
the gap $P_n([u,v))\backslash P_{n+1}([u,v))$ over the interval
$[u,v)$. Since $P_n(I) \supset P_{n+k}(I) \ (k \geq 1)$, there
exists a polynomial $f$ of any order ${\rm deg}(f) \geq 2n - 1$ in
$P_n(I)$.

Finally, by the first statement of the theorem, any polynomial of degree $2n-1$ or $2n$ does not belong to $P_{n+1}(I)$ since $2n-1<2n<2(n+1)-1=2n+1$ and hence if it is in addition a polynomial from
$P_{n}(I)$, then it belongs to the gap
$P_n(I) \setminus P_{n+1}(I)$.
\end{proof}

\section{Truncated Moment Problem and Monotone Matrix Functions.} \label{sec:TrMoPr}

\begin{theorem} \label{the:momentspolyPn}
Let $f(t)=c+b_0 t + b_1 t^2 + \dots + b_{2n-2} t^{2n-1} + b_{2n-1}
t^{2n} + \dots $ be a polynomial of degree at least $2n-1$. Then
\begin{itemize}
\item[{\rm a)}] $M_n(f;0)>0$ if and only if
there is a Borel measure $\mu$ on $\mathbb{R}$ with at least $n$ points in the support, and such that
$$b_k = \int_{\mathbb{R}} t^k d\mu < \infty, \quad \quad \quad (0\leq k \leq 2n-2). $$
Moreover, in this case there exists $\alpha_n > 0$ such that $f \in P_n([0,\alpha_n))$.
\item[{\rm b)}]
If $M_n(f;0) \geq 0$ but $\det M_n(f;0) = 0$, and $r$ is the
smallest positive integer such that $M_{r+1}(f;0)$ is not
invertible, then there exists a Borel measure $\mu$ such that
$$b_k = \int_{\mathbb{R}} t^k d\mu < \infty, \quad \quad \quad (0\leq k \leq 2r-2), $$
and there exists $\alpha>0$ such that $f\in P_{r}([0,\alpha))$.
\end{itemize}
\end{theorem}

\begin{proof}
a) At first we recall that the inequality $M_n(f;0)>0$ means that
the Hankel matrix $M_{n-1}(f;0)$ has a positive Hankel extension
$M_n(f;0)$, and hence by \cite[Theorem 3.9]{CurtoFialkowHJM}, this
is equivalent to the existence of a Borel measure $\mu$ on
$\mathbb{R}$, such that
$$b_k = \int_{\mathbb{R}} t^k d\mu < \infty, \quad \quad \quad (0\leq k \leq 2n-2). $$

Suppose that the measure $\mu$ has at least $n$ points in the support and satisfies $b_k = \int_{\mathbb{R}} t^k d\mu < \infty$ when
$0\leq k \leq 2n-2$. Take arbitrary $n$ points
$t_1, \dots, t_n$ in the support of $\mu$. Then $\mu(I_i)>0$ for any family of $n$ non-overlapping open intervals such that $t_i \in I_i$ for $i=1,\dots,n$. Choose inside each of these open intervals
a closed interval $J_i$ such that
$t_i \in J_i \subset I_i$ and hence also
$\mu (J_i) > 0$ for $i=1,\dots,n$.
For any vector $\vec{c} = (c_0, \dots, c_{n-1}) \in \mathbb{C}^n$, the following holds for the quadratic form
\begin{eqnarray*}
(M_n(f;0)\vec{c} \mid \vec{c}) &=&
\sum_{i=0}^{n-1}\sum_{j=0}^{n-1} b_{i+j} c_j \bar{c_i}  \\
&=&
\sum_{i=0}^{n-1}\sum_{j=0}^{n-1} \int_{\mathbb{R}}
t^{i+j} d \mu \ c_j \bar{c_i}  \\
&=& \int_{\mathbb{R}} \left|
\sum_{i=0}^{n-1} c_i t^{i} \right|^2 d \mu  \\
&\geq & \left|
\sum_{i=0}^{n-1} c_i \hat{t_k}^{i} \right|^2 \mu(J_k) \geq 0,
\end{eqnarray*}
where $\hat{t_k}$ is the minimum point for the continuous function
$\left| \sum_{i=0}^{n-1} c_i t^{i} \right|^2$ on the closed
interval $J_k$. Therefore, the matrix $M_{n}(f;0)$ is at least
positive semi-definite. Moreover, because of $\mu(J_k)>0$, if
$(M_n(f;0)\vec{c} \mid \vec{c})=0$ for some $\vec{c}$, then
$\sum_{i=0}^{n-1} c_i \hat{t_k}^{i}=0$ for all $k=1,\dots, n$.
Since, only the zero polynomial has more roots than its degree,
the only possibility for the linear system to hold is when
$\vec{c}=\vec{0}$. Therefore, the matrix $M_n(f;0)$ is positive
definite ($M_n(f;0)>0$). All elements of this matrix are
polynomials and hence determinants of all submatrices are also
polynomials, and in particular the determinants of all submatrices
with the principal diagonals, consisting from elements of the
principal diagonal of $M_n(f;0)$, are also polynomials and hence
are continuous functions on the real line. There are finitely many
of them and all of them are positive at $t=0$ due to positive
definiteness of $M_n(f;0)$. Each of these polynomials is then
positive on some interval of the form $[0,\alpha)$, and taking the
interval with smallest positive $\alpha$ yields an interval of
this form where the matrix $M_n(f;t)$ is positive. Without
providing a way to construct $[0,\alpha)$ one may alternatively
argue that since eigenvalues of a matrix depend continuously on
each entry of a matrix, there exists a positive $\alpha > 0$ such
that $M_n(f;t)$ is positive semi-definite on $[0,\alpha)$.
Hence $f \in P_n([0,\alpha))$ by \cite[Theorem V, Ch. VIII]{Donoghuebook}.

To prove the converse implication assume that $M_n(f;0)>0$ and let
$\mu$ be a measure satisfying $b_k = \int_{\mathbb{R}} t^k d\mu <
\infty$ when $0\leq k \leq 2n-2$. Assume contrary to the statement
in the theorem, that support of $\mu$ contains less than $n$
points. Let $\{t_1,\dots,t_k\}$, where $k<n$, be the support of
$\mu$. Then there exists a non-zero polynomial $p_n(t)=c_0 + c_1 t
+ \dots + c_{n-1} t^{n-1}$ such that $p_n(t_i)=0$ when $1\leq i
\leq k$. But then for the vector $\vec{c}\neq 0$ of coefficients
of this non-zero polynomial
\begin{eqnarray*}
(M_n(f;0)\vec{c} \mid \vec{c}) = \int_{\mathbb{R}} \left|
\sum_{i=0}^{n-1} c_i t^{i} \right|^2 d \mu = \sum_{j=1}^{k}
\mu(t_j) | p_n(t_j)|^2 =0.
\end{eqnarray*}
This contradicts to the assumption $M_n(f;0)>0$.
Thus $\mu$ must have at least $n$ points in its
support.

b) The existence of the measure such that $b_k = \int_{\mathbb{R}}
t^k d\mu < \infty$ when $0 \leq k \leq 2(r-1)=2r-2$ follows from
\cite[Theorem 3.9]{CurtoFialkowHJM}, and existence of $\alpha>0$
such that $f\in P_{r}([0,\alpha))$ is implied from a) since
$M_{r}(f;0)>0$ by definition of $r$.
\end{proof}

A Hankel rank $rank_h(\vec{\gamma})$ of $\vec{\gamma} = (\gamma_0,
\dots, \gamma_{2k})$ associated to a Hankel matrix
$(\gamma_{i+j})_{i,j=0}^k$ of size $k+1$ is defined as the
smallest integer $i$ obeying $1 \leq i \leq k$ and such that
$\vec{v}_i$ is a linear combination of $\vec{v}_0,\dots,
\vec{v}_{i-1}$, where $\vec{v}_j = (\gamma_{j+l})_{l=0}^{k}$ are
column vectors of the matrix, that is, $(\gamma_{i+j})_{i,j=0}^k =
(\vec{v}_0, \dots ,\vec{v}_k)$. This is a handy notion which we
will use in several examples. According to \cite[Proposition
2.2]{CurtoFialkowHJM} for a positive semidefinite Hankel matrix
the Hankel rank of the defining sequence $\vec{\gamma}$ coincides
with the smallest positive integer $l$ such that the principle
upper left-hand conner submatrix $(\gamma_{i+j})_{i,j=0}^l$ of
size $(l+1)\times (l+1)$ is not invertible, i.e. has zero
determinant, or equivalently this can be rephrased as the largest
integer $l$ such that all the submatrices
$(\gamma_{i+j})_{i,j=0}^{m-1}$ of size $m\times m$ with $1\leq m
\leq l$ are invertible. Thus, in this terminology the integer $r$
we used in the part b) of Theorem \ref{the:momentspolyPn} is
exactly the Hankel rank of the sequence
$\{\gamma_k=\frac{f^{(k+1)}(0)}{(k+1)!}\}_{k=0}^{2(n-1)}$
corresponding to the matrix $M_n(f;0)$.

We present now an example $p$ of a polynomial of  degree $3$ which has determinant of the matrix $M_2(p;t)$ at $t= 0$ being zero, but
$p  \in P_2([0, \alpha))$ for some $\alpha > 0$. 
Let $f(t) = t - t^2 + t^3$. Then
$$
f'(t)  = 1 - 2t + 3t^2, \quad 
f{''}(t) = -2 + 6t, \quad
f{'''}(t) = 6.
$$
Hence
$$
M_2(f;t) = \left(\begin{array}{ll}
1 - 2t + 3t^2&- 1 + 3t\\
-1 + 3t & 1
\end{array}
\right)
$$
Then we have
$$
{\rm det}(M_2(f;t)) = 4t - 6t^2 = -6(t - \frac{1}{3})^2 + \frac{2}{3}.
$$
If take $\alpha = \frac{2}{3}$, $f \in P_2([0, \alpha))$.
\hfill \qed

\vskip 3mm

Next we will show that  a polynomial $p$  of degree $5$  does not
belong to $P_3([0,\alpha))$ for any $\alpha > 0$ when ${\rm
rank}M_3(p;0) = 2$ and $rank_h(\vec{\gamma}) = rank_h(b_0, b_1,
b_2, b_3, b_4) = 1$, where $p(t) = b_0t + b_1t^2 + b_2t^3 + b_3t^4
+ b_4t^5$. Let
\begin{align*}
M_3(p;0) = \left(\begin{array}{ccc}
              b_0&b_1&b_2\\
              b_1&b_2&b_3\\
              b_2&b_3&b_4
              \end{array}
              \right).
\end{align*}
Since $rank_h(\vec{\gamma}) = rank_h(b_0, b_1, b_2, b_3, b_4) =
1$, we have
\begin{align*}
b_1 &= \lambda b_0\\
b_2 &= \lambda b_1 = \lambda^2 b_0\\
b_3 &= \lambda b_2 = \lambda^3 b_0
\end{align*}
for some $\lambda \in \mathbb{R}$. Since ${\rm rank}M_3(p;0) = 2$,
$b_0 > 0$, and we may assume that $b_0 = 1$. Hence we consider
$$
p(t) = t + \lambda t + \lambda^2t^3 + \lambda^3t^4 + ct^5.
$$
for any $c \geq 0$.

\begin{proposition} Let $p(t) = b_0t + b_1t^2 + b_2t^3 + b_3t^4 +
b_4t^5$. Suppose that ${\rm rank}M_3(p;0) = 2$ and
$rank_h(\vec{\gamma}) = rank_h(b_0, b_1, b_2, b_3, b_4) = 1$. Then
there exists no $\alpha
> 0$ that satisfy $p \in P_3([0, \alpha))$.
\end{proposition}

\begin{proof}
{}From the above argument we assume that
$$
p(t) = t + \lambda t + \lambda^2t^3 + \lambda^3t^4 + ct^5.
$$

Consider
$$
M_3(p;t) =
\left(
\begin{array}{ccc}
1 + 2\lambda t + 3\lambda^2t^2 + 4\lambda^3t^3 + 5ct^4
&\lambda + 3\lambda^2 t + 6\lambda^3t^2 + 10ct^3
&\lambda^2 + 4\lambda^3t + 10ct^2\\
\lambda + 3\lambda^2 t + 6\lambda^3t^2 + 10ct^3
&\lambda^2 + 4\lambda^3t + 10ct^2
&\lambda^3 + 5ct\\
\lambda^2 + 4\lambda^3t + 10ct^2
&\lambda^3 + 5ct
&c
\end{array}
\right)
$$
Hence
\begin{align*}
{\rm det}M_3(p;t) &= 30\lambda^4t^2c - 15c^2t^2 - 15\lambda^8t^2 -
30\lambda^5ct^3 + 50\lambda c^2t^3 - 20\lambda^9t^3  \\
& \hspace{4cm} -210\lambda^3c^2t^5 - 105\lambda^6ct^4 - 175c^3t^6\\
&= t^2(30 \lambda^4c - 15c^2 - 15\lambda^8) - 30\lambda^5ct^3  +
50\lambda c^2t^3 - 20\lambda^9t^3 \\ & \hspace{4cm}
-210\lambda^3c^2t^5 - 105\lambda^6ct^4 - 175c^3t^6
\end{align*}

Since
\begin{align*}
30 \lambda^4c - 15c^2 - 15\lambda^8
&= 15(2\lambda^4c - c^2 - \lambda^8)\\
&\leq 15(2\lambda^4c - 2\lambda^4c) = 0,
\end{align*}
where the equality holds when $c = \lambda^4$.

If $c \not= \lambda^4$, the coefficient of $t^2$ is negative, and
there exists $\alpha > 0$, such that $\det M_3(p;t) < 0$ for any $t \in [0,
\alpha)$.

If $c = \lambda^4$, then we have
\begin{align*}
\det M_3(p;t) &= - 210\lambda^3c^2t^5 - 105\lambda^6ct^4
- 175c^3t^6\\
&= - 105\lambda^{10}t^4 - 210\lambda^{11}t^5 - 175\lambda^{12}t^6
\end{align*}
Since $- 105\lambda^{10} < 0$, there exists $\alpha > 0$ such that
$\det M_3(p;t) < 0$ for any $t \in [0, \alpha)$. Hence there
exists no $\alpha > 0$ such that $p \in P_3([0, \alpha))$.
\end{proof}

As another example consider $p(t) = \frac{1}{2}t + t^2 +
\frac{1}{2}t^3 + t^4 + \frac{1}{2}t^5$. Then the matrix ${\rm
det}M_3(p;0) = 0$ and $M_3(p;0)$ has rank $2$. Note that ${\rm
rank}_h(\frac{1}{2},1,\frac{1}{2},1,\frac{1}{2}) = 2$. Therefore,
the situation is different in the previous proposition. Since
$$
M_3(p;t) = \left(
\begin{array}{ccc}
\frac{1}{2} + 2t + \frac{3}{2}t^2 + 2t^3 + \frac{5}{2}t^4 &1 +
\frac{3}{2}t + 6t^2 + 5t^3
&\frac{1}{2} + 4t + 5t^2\\
1 + \frac{3}{2}t + 6t^2 + 5t^3 &\frac{1}{2} + 4t + 5t^2
&1 + \frac{5}{2}t\\
\frac{1}{2} + 4t + 5t^2 &1 + \frac{5}{2}t &\frac{1}{2}
\end{array}
\right),
$$
we have
\begin{align*}
{\rm det}M_3(p;t) = \frac{9}{2}t + \frac{63}{8}t^2 -
\frac{27}{2}t^3 - \frac{93}{2}t^4 - 45t^5 - \frac{175}{8}t^6.
\end{align*}
Hence there exists $\alpha > 0$ such that $p \in P_3([0,\alpha))$.

The previous results and examples imply the following theorem
which is concerned with catching conditions for a more precise
determination of position of a given function with respect to the
decreasing sequences of inclusions for the classes of matrix
monotone functions.

\begin{theorem}
Let $0\in [0,\alpha)$ and let $f$ be a polynomial such that $f\in
P_n([0,\alpha))$.
\begin{enumerate}
\item If $f\in P_{n+1}([0,\alpha))$,
then there exists a Borel measure $\mu$ such that $b_k =
\int_{\mathbb{R}} t^k d \mu $ for $ 0\leq k \leq 2n-1$.
\item If $f\in P_{n+1}([0,\alpha))$ and $M_{n+1}(f;0)>0$, then there exists Borel measure $\mu$ such that
$b_k = \int t^k d \mu $ for $ 0\leq k \leq 2n$;
\item
Let $f\in P_n([0,\alpha))$ and let $r$ be the smallest number such
that the submatrix $M_{r+1}(f;0)$ is not invertible. If $r < {\rm
rank} (M_n(f;0))$, then $f \not \in P_{n+1}([0,\alpha)).$
\end{enumerate}
\end{theorem}
\begin{proof}
1) If $f\in P_{n+1}([0,\alpha))$, then $M_{n+1}(f;0)\geq 0$. Hence
$M_{n}(f;0)$ has a positive semidefinite Hankel extension, and
thus according to  \cite[Theorem 3.1]{CurtoFialkowHJM} there
exists a Borel measure $\mu$ such that $b_k = \int_{\mathbb{R}}
t^k d \mu $ for $ 0\leq k \leq 2n-1$.

2) If $f\in P_{n+1}([0,\alpha))$ and moreover $M_{n+1}(f;0)> 0$, then the
existence of a Borel measure $\mu$ such that $b_k =
\int_{\mathbb{R}} t^k d \mu $ for $ 0\leq k \leq 2n-1$ is already
secured by 1) and then the fact that the next coefficient $b_{2n}$
may also be determined by the moment $b_{2n}=\int_{\mathbb{R}}
t^{2n} d \mu  $, as claimed in the theorem, follows from the
statement a) of Theorem 3, since in this case $2(n+1)-2=2n$.

3) If $f\in P_{n+1}([0,\alpha))$, then $M_{n}(f;0)$ has positive
Hankel extension, and by \cite[Theorem 3.9]{CurtoFialkowHJM} the
ordinary matrix rank of $M_n(f;0)$ has to be  equal to the Hankel
rank $r$. Hence, if this equality does not hold, then $f \not \in
P_{n+1}([0,\alpha))$.
\end{proof}

\section{Rank and the Hadamard product.}

In this section we treat the rank comparison problem
between an $n \times n$ matrix $A$
and $A \bullet D$, where $\bullet$ means the Hadmard product of the matrix $A$ with another matrix $D$,
and then we show that these matrix results are useful for
understanding of the transformations of matrix monotonicity
properties of functions when changing from one interval to another.

Let $A$ be an $n \times n$ matrix
$$
\left(\begin{array}{ccc}
a_{11}&\cdots&a_{1n}\\
\vdots&&\vdots\\
a_{n1}&\cdots&a_{nn}
\end{array}
\right).
$$
By $A(k)$ for $(1 \leq k \leq n)$ we denote the
$k \times k$ left upper conner matrix, that is,
$$
A(k) =
\left(\begin{array}{ccc}
a_{11}&\cdots&a_{1k}\\
\vdots&&\vdots\\
a_{k1}&\cdots&a_{kk}
\end{array}
\right).
$$

For $n \times n$ matrices $A = (a_{ij})$ and $D = (d_{ij})$ we
write
$$
A \bullet D = (a_{ij}d_{ij}).
$$

\begin{lemma} \label{lemma:rankineqHadPr}
(\cite[Theorem 5.1.7]{HJ})
Let $A$ and $B$ be two $n \times n$ matrices.
Then
$$
{\rm rank}(A \bullet B) \leq ({\rm rank}(A))({\rm rank}(B)).
$$
\end{lemma}

For an $n \times n$ matrix $A$ with real eigenvalues,
we denote by $\lambda_{min}(A)$ the minimal
eigenvalue of $A$, and by
$\lambda_{max}(A)$ the maximal eigenvalue of $A$.

\begin{lemma} \label{lemma2:eigHadProd} (\cite[Theorem 5.3.4]{HJ})
Let
$A$ and $B$ be two $n \times n$ positive semidefinite matrices.
Then any eigenvalue $\lambda(A \bullet B)$ of $A \bullet B$ satisfies
\begin{align*}
\lambda_{min}(A)\lambda_{min}(B)
&\leq (\min_{1 \leq i \leq n}(a_{ii}))\lambda_{min}(B)\\
&\leq \lambda(A \bullet B)\\
&\leq (\max_{1 \leq i \leq n}(a_{ii}))\lambda_{max}(B)\\
&\leq \lambda_{max}(A)\lambda_{max}(B)
\end{align*}
\end{lemma}

Using the above two Lemmas we will show the following result.

\begin{proposition}
Let $A$ and $D$ be positive semidefinite $n \times n$ matrices.
Suppose that $D = [\alpha^{i+j-1}]$ and $\alpha > 0$.
Then for $1 \leq k \leq n$, the matrix
$A(k)$ is invertible if and only if $(A \bullet D)(k)$ is invertible.
\end{proposition}

\begin{proof}
Note that ${\rm rank}(D(k)) = 1$ for all
$1 \leq k \leq n$.

Suppose that $A(k)$ is invertible. Then since $(A \bullet D)(k) =
A(k) \bullet D(k) = D(k) \bullet A(k)$, from Lemma
\ref{lemma2:eigHadProd},
$$
\lambda((A \bullet D)(k)) \geq (\min_{1 \leq i \leq n}\{d_{ii}\})
\lambda_{min}(A(k)) > 0.
$$
Since any eigenvalue of  $(A \bullet D)(k)$ is positive,
$(A \bullet D)(k)$ is invertible.

Conversely, suppose that $(A \bullet D)(k)$ is invertible.
From Lemma \ref{lemma:rankineqHadPr}
\begin{align*}
k = {\rm rank}((A \bullet D)(k)) &\leq {\rm rank}(A(k)){\rm rank}(D(k))\\
&\leq {\rm rank}(A(k)).
\end{align*}

Hence ${\rm rank}(A(k)) = k$, and $A(k)$ is invertible.
\end{proof}

\begin{corollary} \label{cor:HPr1}
Let $A$ be a positive semidefinite $n \times n$ matrix
and $D = (\alpha^{i+j-1})_{i,j=1}^{n}$ with  $\alpha \not= 0$ be a
$n \times n$ matrix.
Then for $1 \leq k \leq n$,
$A(k)$ is invertible if and only if $(A \bullet D)(k)$ is invertible.
\end{corollary}

\begin{proof}
If $\alpha > 0$, then $D$ is positive semidefinite matrix from the elementary
calculation.
So, the conclusion follows from the previous proposition.

If $\alpha < 0$, $- D$ is positive semidefinite.
Since $- (A \bullet D) = A \bullet (- D)$ and
${\rm rank}(A \bullet D) = {\rm rank}(- (A \bullet D))$, we get the
conclusion.
\end{proof}

The presented results on the rank for Hadamard product of matrices are quite useful when
attempting to describe
how the classes $P_n(I)$ are
related to each other for different intervals.

Let $f(t) = b_0 t + b_1 t^2 + \dots + b_{2n-2} t^{2n-1}$. The
interval $[0,a)$ is transformed bijectively to $[u,v)$ by the
operator monotone affine mapping  $h(t)=\frac{v-u}{a}t +u$ with
the operator monotone composition inverse $h^{\circ (-1)}(t) =
\frac{a}{v-u}t-\frac{au}{v-u}$. Then applying the chain rule and
the affine form of $h$ and $h^{\circ (-1)}$ we have $M_n(f\circ
h^{\circ (-1)};0) = M_n(f;a)\bullet
\left((\frac{a}{v-u})^{i+j-1})_{i,j=0}^{n}\right)$ and hence the
rank of $M_n(f\circ h^{\circ (-1)};0)$ and of $M_n(f;a)$ coincide
according to Corollary \ref{cor:HPr1}. Thus if $u(t)=c_0 (t-a)
+\dots +c_{2n-1}(t-a)^{2n-1}$ on the interval $[u,v)$, and
correspondingly $u\circ h^{\circ (-1)}(t) = b_0 t + \dots +
b_{2n-2} t^{2n-1}$ on the interval $[0,\alpha)$, then there exists
a measure $\mu$ such that $b_k = \int t^k d \mu $ for $ 0\leq k
\leq 2n-2$ if and only if there exists a measure $\tilde{\mu}$
such that $c_k = \int t^k d \tilde{\mu} $ for $ 0\leq k \leq
2n-2$. Therefore, in this sense there is a correspondence between
the structure of those polynomials in $P_n([u,v))$ and
$P_n([0,\alpha))$.

\section{The characterization of operator monotone functions over $[0,a)$}


Let $I = [0, a)$ for $a > 0$.

\begin{definition} \label{def:Mnklas} Let ${\cal M}_n(I)$ be the class of functions such that
$f \in {\cal M}_n(I)$ if for all $a_k \in \mathbb{R}$, $\lambda_k
\in (0, a)$ for $1 \leq k \leq 2n$
$$
\left\{\begin{array}{ll}
\sum_{k=1}^{2n}a_k\frac{\lambda_k}{t + \lambda_k} \geq 0 \
\hbox{for} \ t > 0\\
\sum_{k=1}^{2n}a_k = 0\\
\end{array}
\right.
$$
implies that
$$
\sum_{k=1}^{2n}a_kf(\lambda_k) \geq 0.
$$
\end{definition}
The above class ${\cal M}_n(I)$ is a finite interval version of
the class ${\cal M}_n$ in \cite{SparrMathScandart}.

\begin{remark} \label{rem:Mnklascond}
Since for $t > 0$
\begin{align*}
\sum_{k=1}^{2n}a_k\frac{\lambda_kt - 1}{t + \lambda_k} =
(t + \frac{1}{t})\sum_{k=1}^{2n}a_k\frac{\lambda_k}{t + \lambda_k},
\end{align*}
$f \in {\cal M}_n(I)$ if and only if
$$
\left\{\begin{array}{ll}
\sum_{k=1}^{2n}a_k\frac{\lambda_kt - 1}{t + \lambda_k} \geq 0 \
\hbox{for} \ t > 0\\
\sum_{k=1}^{2n}a_k = 0\\
\end{array}
\right.
$$
implies that
$$
\sum_{k=1}^{2n}a_kf(\lambda_k) \geq 0.
$$
\hfill\qed
\end{remark}

We use the following characterization of monotone functions in
$P_n(I)$ as in \cite{SparrMathScandart}.

\begin{lemma}
\label{lemmaPnkvadformkrit} For $\alpha$ and $x$ in $\mathbb{C}^n$
set $\|x\|_\alpha = \left(\sum_{k=1}^n\alpha_k|x_k|^2
\right)^\frac{1}{2}$.

Then
$f \in P_n(I)$ if and only if
\begin{align*}
\hbox{for all } n \times n \ \hbox{unitary} \ U \ \hbox{with}\
&\|U\|_{\alpha,\beta} \leq 1 \
(\alpha, \beta \in \mathbb{C}^n \cap I^n):\\
& \sum_{k=1}^n f(\alpha_k)|x_k|^2 \geq \sum_{k=1}^n
f(\beta_k)|(Ux)_k|^2, \ \forall x \in \mathbb{C}^n,
\end{align*}
where
$$
\|U\|_{\alpha,\beta} = \sup_{\begin{array}{c} x \in
\mathbb{C}^n \\[-3mm]
x \not= 0
\end{array}}
\frac{\|Ux\|_\beta}{\|x\|_\alpha}.
$$
\end{lemma}

\begin{proof}
Let $A$ and $B$ be two hermitian $n \times n$ matrices with
eigenvalues contained in $I$. Then
\begin{align*}
A \geq B &\Leftrightarrow \sum_{k=1}^n\alpha_k|x_k|^2 \geq
\sum_{k=1}^n\beta_k|(Ux)_k|^2,  \
\forall x \in \mathbb{C}^n
\end{align*}
where $x = \left(x_1, x_2, \hdots, x_n \right)^T $, $\alpha_1,
\alpha_2, \dots, \alpha_n$ are eigenvalues for $A$, $\beta_1,
\beta_2, \dots, \beta_n$ eigenvelues for $B$, and $U$ is an
appropriate  $n \times n$  unitary. Every unitary arises for some
choice of $A$ and $B$. Hence we have
\begin{align*}
f \in P_n(I) &\Leftrightarrow
\forall \  n \times n \ \hbox{unitary} \ U \ \hbox{with}\
\|U\|_{\alpha,\beta} \leq 1 \
(\alpha, \beta \in \mathbb{C}^n \cap I^n)\\
& \sum_{k=1}^n f(\alpha_k)|x_k|^2 \geq \sum_{k=1}^n
f(\beta_k)|(Ux)_k|^2, \ \forall x \in \mathbb{C}^n.
\end{align*}
\end{proof}

As for ${\cal M}_n$ in \cite{SparrMathScandart} we have the
following fundamental inclusion.

\begin{proposition}
For all $n \in \mathbb{N}$
$$
P_{n+1}(I) \subseteq {\cal M}_n(I) \subseteq P_n(I).
$$
\end{proposition}

\begin{proof} Using the same argument in \cite{SparrMathScandart}
we can show the inclusion ${\cal M}_n(I) \subset P_n(I)$.

To show the inclusion $P_{n+1} \subset {\cal M}_n(I)$ we take the
same steps in \cite{SparrMathScandart}. Take an arbitrary $f \in
P_{n+1}(I)$, and choose $0 < \lambda_1 < \lambda_2 < \cdots <
\lambda_{2n} < a$. Consider
$$
\frac{p(t)}{\pi(t)} = \sum_{k=1}^{2n}a_k \frac{\lambda_k}{t +
\lambda_k},
$$
where $p$ is any polynomial of degree less than or equal to $2n - 1$
(write the class of such polynomials  by $Pol(2n - 1)$)
such that $p(t) \geq 0$ for $t > 0$ and $p(0) = 0$, and
$$
\pi(t) = \Pi_{i=1}^{2n}(t + \lambda_i).
$$
Then we may show that
\begin{equation} \label{proofMnPnineq1}
\sum_{k=1}^{2n}a_kf(\lambda_k) \geq 0.
\end{equation}
Note that
$$
a_k = \frac{p(-\lambda_k)}{\lambda_k\pi'(-\lambda_k)}
$$
for $1 \leq k \leq 2n$, and the polynomials with the above
property can be written as
$$
p(t) = tq_1(t)^2 + q_2(t)^2,
$$
where $q_1, q_2 \in Pol(n - 1)$ and $q_2(0) = 0$. (For example see
\cite[Lemma 7.6.1]{Pe}.) Hence, because of linearity we only have
to consider the two cases $p(t) = tq(t)^2$ and $p(t) = q(t)^2$
with $q(0) = 0$.

When $p(t) = tq(t)^2$, we can show the inequality
\eqref{proofMnPnineq1} by the same argument as in (i) of the proof
in \cite[Lemma 1]{SparrMathScandart}. We write $0 < \lambda_1 <
\lambda_2 < \cdots < \lambda_{2n} < a$ as $0 < \beta_1 < \alpha_1
< \beta_2 < \cdots < \beta_n < \alpha_n < a$. When $p(t) =
q(t)^2$, $q \in Pol(n - 1)$, $q(0) = 0$, write
$$
\frac{q^2(t)}{\pi(t)}
= \sum_{k=1}^ny_k^2\frac{\beta_k}{t + \beta_k} -
\sum_{k=1}^nx_k^2\frac{\alpha_k}{t + \alpha_k},
$$
where
$$
y_k^2 = \frac{q^2(-\beta_k)}{\beta_k\pi'(-\beta_k)}, \
x_k^2 = \frac{- q^2(-\alpha_k)}{\alpha_k\pi'(-\alpha_k)}.
$$
We extend $0 < \beta_1 < \alpha_1 < \beta_2 < \cdots < \beta_n <
\alpha_n < a$ with $\delta$ and $\omega$ such that
$$
0 < \delta < \beta_1 < \alpha_1 < \beta_2 < \cdots
< \beta_n < \alpha_n < \omega < a
$$
 and consider
 $$
\frac{t}{t + \delta} \cdot \frac{q^2(t)}{\pi(t)} \cdot
\frac{t+a}{t+\omega}.
$$
Note that
$$
\frac{t}{t + \delta} \cdot \frac{q^2(t)}{\pi(t)} \cdot
\frac{t+a}{t+\omega} \ \rightarrow \frac{q^2(t)}{\pi(t)}
$$
as $\delta \rightarrow 0$ and $\omega \rightarrow a$. By the
partial fraction expansion
\begin{equation} \label{proofMnPnineq2}
\frac{t}{t + \delta} \cdot \frac{q^2(t)}{\pi(t)} \cdot
\frac{t+a}{t+\omega} = - \tilde{x}_0^2\frac{\delta}{t+\delta} -
\sum_{k=1}^n\tilde{x}_k^2\frac{\alpha_k}{t+\alpha_k} +
\sum_{k=1}^n\tilde{y}_k^2\frac{\beta_k}{t+\beta_k} +
\tilde{y}_{n+1}^2\frac{\omega}{t+\omega},
\end{equation}
where $\tilde{x}_k$ and $\tilde{y}_k$ are defined
similar as $x_k$ and $y_k$, and
$$
\tilde{x}_k \rightarrow x_k, \ \tilde{y}_k \rightarrow y_k, \ \ \
1 \leq k \leq n
$$
as $\delta \rightarrow 0$ and $\omega \rightarrow a$.
Moreover
$$
\tilde{x}_0^2 = \frac{q^2(-\delta)}{\pi(-\delta)} \cdot \frac{a -
\delta}{\omega - \delta} = {\mathcal O}(\delta^2), \delta
\rightarrow 0,
$$
since $q(0) = 0$, and
$$
\tilde{y}_{n+1}^2 = \frac{q^2(-\omega)}{\pi(-\omega)} \cdot
\frac{-\omega+a}{\omega-\delta} = {\mathcal O}(-\omega+a), \
\omega \rightarrow a.
$$
Let $f \in P_{n+1}(I)$. By letting $t = 0$ in
\eqref{proofMnPnineq2} we have
\begin{eqnarray}
\label{proofMnPnineq3} - \tilde{x}_0^2 - \sum_{k=1}^n\tilde{x}_k^2
+ \sum_{k=1}^n\tilde{y}_k^2 + \tilde{y}_{n+1}^2 = 0.
\end{eqnarray}
Since $\frac{t}{t +
\delta}\frac{q^2(t)}{\pi(t)}\frac{t+a}{t+\omega}t \geq 0$,
$$
\lim_{t\rightarrow\infty} \frac{t}{t + \delta} \cdot
\frac{q^2(t)}{\pi(t)} \cdot \frac{t+\alpha}{t+\omega} \cdot t \geq
0,
$$
Hence
$$
-\tilde{x}_0^2\delta - \sum_{k=1}^n\tilde{x}_k^2\alpha_i
+ \sum_{k=1}^n\tilde{y}_k^2\beta_k + \tilde{y}_{n+1}^2\omega \geq 0.
$$
Since $f \in P_{n+1}(I)$ and \eqref{proofMnPnineq3},
$$
-\tilde{x}_0^2f(\delta) - \sum_{k=1}^n\tilde{x}_k^2f(\alpha_k)
+ \sum_{k=1}^n\tilde{y}_k^2f(\beta_k) + \tilde{y}_{n+1}^2f(\omega) \geq 0.
$$
(See Lemma \ref{lemmaPnkvadformkrit} and $(1')$ in
\cite{SparrMathScandart}.)

By the same argument as in (ii) in the proof in \cite[Lemma
1]{SparrMathScandart}, we have
$$
\lim_{\delta\rightarrow 0}\tilde{x}_0^2f(\delta)
= \lim_{\omega\rightarrow a}\tilde{y}_{n+1}^2f(\omega)
= 0.
$$
Both equalities come from the same proof as in \cite[Lemma
1]{SparrMathScandart}. Indeed, we consider the following
inequality, which is  used in \cite[Lemma 1]{SparrMathScandart}.
\begin{align*}
&- \frac{(c - \beta_1)^2}
{(\alpha_1 - \beta_1)(\beta_2 - \beta_1)(\alpha_2 - \beta_1)}
f(\beta_1)
- \frac{(c - \beta_2)^2}
{(\beta_1 - \beta_2)(\alpha_1 - \beta_2)(\alpha_2 - \beta_2)}
f(\beta_2)
\\
&+ \frac{(c - \alpha_1)^2}
{(\beta_1 - \alpha_1)(\beta_2 - \alpha_1)(\alpha_2 - \alpha_1)}
f(\alpha_1)
+ \frac{(c - \alpha_2)^2}
{(\beta_1 - \alpha_2)(\alpha_1 - \alpha_2)(\beta_2 - \alpha_2)}
f(\alpha_2) \geq 0.
\end{align*}
This comes from the fact that $f \in P_2(I)$ and Lemma \ref{lemmaPnkvadformkrit}.
(See $(1')$ in \cite{SparrMathScandart}.)

To get the first equality, set $c = \beta_1 = \frac{\delta}{2}$,
$\alpha_1 = \delta, \beta_1 = \frac{a}{4}$, and $\alpha_2 = \frac{a}{2}$.
Then we have
\begin{align*}
\delta f(\delta)
\geq
&\frac{8}{a}\left\{(\frac{a}{2} - \delta)(\frac{a}{4} - \frac{\delta}{2})
f(\frac{a}{4})
- (\frac{a}{2} - \frac{\delta}{2})(\frac{a}{4} - \delta)f(\frac{a}{2})
\right\}
\end{align*}
Hence
\begin{align*}
\lim\inf_{\delta\rightarrow 0}\delta^2f(\delta) \geq 0.
\end{align*}
Since $f$ is monotone, $\delta^2f(\delta) \leq \delta^2f(\frac{a}{2})$
for $\delta < \frac{a}{2}$. Then

$$
\lim\sup_{\delta\rightarrow 0}\delta^2f(\delta) \leq
\lim\sup_{\delta \rightarrow 0}\delta^2 f(\frac{a}{2}) = 0,
$$
hence
$\lim_{\delta\rightarrow 0}\delta^2f(\delta) = 0$.
Therefore we have
\begin{align*}
\lim_{\delta\rightarrow 0}\tilde{x}_0^2f(\delta)
= \lim_{\delta\rightarrow 0}\frac{\tilde{x}_0^2}{\delta^2}\delta^2f(\delta) = 0
\ (\tilde{x}_0^2 = {\mathcal O}(\delta^2)) \end{align*}

To get the second equality
set $c = \alpha_1 = \frac{a}{4}$, $\beta_1 = \frac{a}{8}$,
$\beta_2 = \frac{\omega}{2}$, and $\alpha_2 = \omega$
with $0 < \beta_1 < \alpha_1 < \beta_2 < \alpha_2 < a$.
Then we have
\begin{align*}
f(\omega) & \leq -\frac{\frac{a}{8}(\omega -
\frac{a}{8})\frac{\omega}{2}} {\omega - \frac{a}{4}}f(\frac{a}{8})
- (\frac{\frac{\omega}{2} - \frac{a}{4}) (\omega -
\frac{a}{8})\frac{\omega}{2}} {(\omega -
\frac{a}{4})(\frac{\omega}{2} - \frac{a}{8})\frac{\omega}{2}}
f(\frac{\omega}{2})\   (= h(\omega))
\end{align*}
for $\omega \in (\frac{a}{2}, a)$.

Multiplying the above inequality by $(-\omega + a) > 0$
for $\omega \in (0,a)$ we obtain that
$$
 (-\omega + a)f(\omega) \leq (-\omega + a)h(\omega) \ \ \omega \in (0, a).
$$
Hence
$$
\lim_{\omega\rightarrow a}\sup(a - \omega)f(\omega) \leq 0.
$$

On the contrary, since $(a - \omega)f(\omega) \geq (a - \omega)f(\frac{a}{2})$
for $\omega \in (\frac{a}{2}, a)$,
$$
\lim_{\omega\rightarrow a}\inf(a - \omega)f(\omega) \geq 0,
$$
and hence $\lim_{\omega\rightarrow a} (a-w)f(\omega) = 0$.
Therefore, we have
\begin{align*}
\lim_{\omega\rightarrow a}\tilde{y}_{n+1}^2
=
\lim_{\omega\rightarrow a}\frac{\tilde{y}_{n+1}^2}{(a - \omega)}
(a - \omega)f(\omega) = 0.\ (\tilde{y}_{n+1}^2 = {\mathcal O}(a - \omega))
\end{align*}
Hence, we get
$$
- \sum_{k=1}^nx_k^2f(\alpha_k) + \sum_{k=1}^ny_k^2f(\beta_k) \geq 0,
$$
and $f \in {\cal M}_n(I)$.
\end{proof}

From the above inclusion property, we have the following characterization of operator monotone functions. 
\begin{theorem}
The function $f$ is operator monotone on $I$ if and only if
$$
f \in \cap_{n=1}^\infty {\cal M}_n(I).
$$
\end{theorem}

\subsection{Examples}

Let $g_n$ be polynomials considered in
\cite{HansenJiTomiyama-art}. We show that $g_n \in
P_n([0,\alpha_n])\backslash {\cal M}_n([0,\alpha_n])$ for some
$\alpha_n > 0$ and $n = 2, 3, 4, 5$ using Maple.

We believe that $g_n \in P_n([0,\alpha_n])\backslash {\cal
M}_n([0,\alpha_n])$ for some $\alpha_n > 0$ and arbitrary $n \geq
2$.

\subsubsection{$g_2$ case}
Let $g_2(x) = x + 1/3x^3$ and let $M_2(g_2;x)$ be the matrix
function corresponding to $g_2$, that is,
$$
M_2(g_2;x) = \left(\begin{array}{cc}
1+x^2&x\\
x&1/3\\
\end{array}
\right)
$$

We claim that $g_2 \in P_2([0,\alpha_2])$ for some $\alpha_2 >
\frac{1}{2}$. To this end we have only to show that $M_2(g_2;x)$
is positive definite for all $x \in [0,\frac{1}{2}]$. Since the
determinant $\det(M_2(g_2;x))$ is $1/3 - 2/3x^2$, it is easily
seen that $\det(M_2(g_2;x)) > 0$ for all $x \in [0, 1/2]$. Hence
$\alpha_2 > 1/2$ and $g_2 \in P([0,\alpha_2])$ by \cite[Theorem
VIII.V]{Donoghuebook}.

Next we show that $g_2 \notin {\cal M}([0,\alpha_2])$. We take the
polynomial $p(x) = x^2$ and
$$
\lambda_k = k/8, \ \ 1 \leq k \leq 4.
$$
Since $a_k = \frac{p(-\lambda_k)}{\lambda_k\pi'(-\lambda_k)}$
for $1 \leq k \leq 4$ where
$\pi(x) = \Pi_{k=1}^4(x + \lambda_k)$, we have
$$
\sum_{k=1}^4a_kg_2(\lambda_k) = -\frac{5}{12} < 0.
$$
This implies that $g_2 \notin M_2([0,\alpha_2])$ by Definition \ref{def:Mnklas}.
\hfill\qed

\subsubsection{$g_3$ case}

Let $g_3(x) = x + 1/3x^3 + 1/5x^5$ and let $M_3(g_3;x)$ be the
corresponding matrix function for $g_3$, that is,
$$
M_3(g_3;x) = \left(\begin{array}{ccc}
1+x^2+x^4&x+2x^3&1/3+2x^2\\
x+2x^3&1/3+2x^2&x\\
1/3+2x^2&x&1/5\\
\end{array}
\right).
$$
We claim that $g_3 \in P_3([0,\alpha_3])$ for some $\alpha_3 >
\frac{1}{5}$. To get this we have only to show that $M_3(g_3;x)$
is positive definite for all $x \in [0,\frac{1}{5}]$. The
determinants of principal matrices of $M_3(g_3;x)$ are as follows:
\begin{align*}
&(3.1) \det(M_3(g_3;x)_{22}) = 1/3 + 4/3x^2 - 5/3x^4 - 2x^6\\
&(3.2) \det(M_3(g_3;x)) = 4/135 - 11/15x^2 - 7/5x^6,
\end{align*}
where $M_3(g_3;x)_{22}$ means the $2 \times 2$upper part of
$M_3(g_3;x)$. Then we can conclude that
\begin{align*}
&\det(M_3(g_3;x)_{22}) > 0\\
&\det(M_3(g_3;x)) > 0
\end{align*}
for all $x \in [0, \frac{1}{5}]$ from graphs in Appendix.

From the above two graphs we can conclude that $M_3(g_3;x)$ is
positive definite for any $x \in [0, \frac{1}{5}]$. (See
\cite[Theorem I.3.3]{Donoghuebook} for example.) Hence $g_3 \in
P_3([0,\alpha_3])$ for some $\alpha_3 > \frac{1}{5}$.

Next we show that $g_3 \notin {\cal M}_3([0,\alpha_3])$. We take
the polynomial $p(x) = x^4$ and
$$
\lambda_k = k/30, \ \ 1 \leq k \leq 6.
$$
Since $a_k = \frac{p(-\lambda_k)}{\lambda_k\pi'(-\lambda_k)}$
for $1 \leq k \leq 6$ and
$\pi(x) = \Pi_{k=1}^6(x + \lambda_k)$, we have
$$
\sum_{k=1}^6a_kg_2(\lambda_k) = -\frac{1897}{7500} < 0.
$$
This implies that $g_3 \notin {\cal M}_3([0,\alpha_3])$ by
Definition \ref{def:Mnklas}. \hfill\qed

\subsubsection{$g_4$ case}

Let $g_4(x) = x + 1/3x^3 + 1/5x^5 + 1/7x^7$ and let $M_4(g_4;x)$
be the corresponding matrix function for $g_4$, that is,
$$
M_4(g_4;x) = \left(\begin{array}{cccc}
1+x^2+x^4+x^6&x+2x^3+3x^5&1/3+2x^2+5x^4&x+5x^3\\
x+2x^3+3x^5&1/3+2x^2+5x^4&x+5x^3&1/5+3x^2\\
1/3+2x^2+5x^4&x+5x^3&1/5+3x^2&x\\
x+5x^3&1/5+3x^2&x&1/7
\end{array}
\right).
$$
We claim that $g_4 \in P_4([0,\alpha_4])$ for some $\alpha_4 >
\frac{1}{25}$. To this end we have only to show that $M_4(g_4;x)$
are positive definite for all $x \in [0,\frac{1}{25}]$. The
determinants of all principal matrices of $M_4(g_4;x)$ are as
follows:
\begin{align*}
&(4.1)\ \det(M_4(g_4;x)_{11}) = 1+x^2+x^4+x^6\\
&(4.2)\ \det(M_4(g_4;x)_{22}) = 1/3+4/3x^2+10/3x^4-8/3x^6-5x^8-4x^10\\
&(4.3)\ \det(M_4(g_4;x)_{33}) = 4/135+4/15x^2-10/3x^4-118/15x^6+2x^8-54/5x^{10}-12x^{12}.\\
&(4.4)\ \det(M_4(g_4;x)) = -848/7875x^2 + 72/7x^{12} - 188/175x^8
+
72/35x^{10} + 16/23625 + \\
& \hspace{7.3cm} +1472/7875x^4 - 4712/875x^6.
\end{align*}
Hence we can conclude that  $M_4(g_4;x)$ are positive  definite
for all $x \in [0,\frac{1}{25}]$ using the Maple, because each of
the determinants is strictly positive for any $x \in [0,
\frac{1}{25}]$. (See graphs in the Appendix.)

Next we claim that $g_4 \not\in {\cal M}_4([0,\alpha_4])$. We take
$p(x) = x^6$ and $\lambda_k = \frac{k}{200}$ for $1 \leq k \leq
8$. Then by the same argument as in the case of $n = 2$ and $n =
3$ we have
$$
\sum_{k=1}^8a_kg_4(\lambda_k) = - \frac{33766394903}{56\times 10^{10}} < 0,
$$
hence $g_4 \not\in {\cal M}_4([0,\alpha_4])$. \hfill\qed


\subsubsection{$g_5$ case}

Let $g_5(x) = x + 1/3x^3 + 1/5x^5 + 1/7x^7 + 1/9x^9$ and let
$M_5(g_5;x)$ be the corresponding matrix function for $g_5$, that
is,
\begin{align*}
& M_5(g_5;x)
=\\
& {\tiny
\left(\begin{array}{ccccc}
1+x^2+x^4+x^6+x^8&x+2x^3+3x^5+4x^7&1/3+2x^2+5x^4+28/3x^6&x+5x^3+14x^5&1/5+3x^2+14x^4\\
x+2x^3+3x^5+4x^7&1/3+2x^2+5x^4+28/3x^6&x+5x^3+14x^5&1/5+3x^2+14x^4&x+28/3x^3\\
1/3+2x^2+5x^4+28/3x^6&x+5x^3+14x^5&1/5+3x^2+14x^4&x+28/3x^3&1/7+4x^2\\
x+5x^3+14x^5&1/5+3x^2+14x^4&x+28/3x^3&1/7+4x^2&x\\
1/5+3x^2+14x^4&x+28/3x^3&1/7+4x^2&x&1/9\\
\end{array}
\right). }
\end{align*}
We claim that $g_5 \in P_3([0,\alpha_5])$ for some $\alpha_5 >
0.032$. To this end we have only to show that all principal
matrices of $M_5(g_5;x)$ are positive definite for all $x \in
[0,0.032]$. The determinants of principal  matrices of
$M_5(g_5;x)$ are as follows:
\begin{align*}
&(5.1)\ \det(M_5(g_5;x)_{11}) = 1+x^2+x^4+x^6+x^8\\
&(5.2)\ \det(M_5(g_5;x)_{22}) =
\frac{1}{3}+\frac{4}{3}x^2+\frac{10}{3}x^4
+\frac{20}{3}x^6-\frac{10}{3}x^8-\frac{26}{3}x^{10}
-\frac{29}{3}x^{12}-\frac{20}{3}x^{14}\\
&(5.3)\ \det(M_5(g_5;x)_{33})\\
& = \frac{4}{15}x^2+\frac{4}{3}x^4-\frac{82}{9}x^6
-\frac{97}{3}x^8-\frac{656}{15}x^{10}+\frac{613}{45}x^{12}
-42x^{14}-\frac{242}{3}x^{16}-\frac{1540}{27}x^{18}+\frac{4}{135}\\
&(5.4)\ \det(M_5(g_5;x)_{44}) \\
&= \frac{256}{23625}x^2-\frac{18824}{23625}x^4-\frac{7136}{2625}x^6+
\frac{5588}{875}x^8+\frac{16}{23625}
+\frac{6776}{27}x^{20}-\frac{137576}{1575}x^{10}\\
&-\frac{254962}{1575}x^{12}+\frac{44}{3}x^{14}-
\frac{2728}{945}x^{16}+\frac{6776}{27}x^{18}\\
&(5.5) \det(M_5(g_5;x)) \\
&= - \frac{34256}{10418625}x^2 + \frac{69212}{243}x^{20}
+ \frac{20251814}{138915}x^{12} + \frac{1024}{260465625}
+ \frac{1424236}{694575}x^8 - \frac{284372}{138915}x^{10}\\
&+ \frac{216592}{10418625}x^4 + \frac{173030}{1701}x^{18}
- \frac{644930}{11907}x^{16} - \frac{1213916}{694575}x^6 +
\frac{1617407}{27783}x^{14}.
\end{align*}

Hence we can conclude that  $M_5(g_5;x)$ are positive  definite
for all $x \in [0,\frac{1}{125}]$ using Maple, because each
determinants is strictly positive for any $x \in [0, 0.032]$. (See
graphs in the Appendix.)

Next we claim that $g_5 \not\in {\cal M}_5([0,\alpha_5])$. As in
the case of $g_5$, we shall find a polynomial $p$, and positive
number $\lambda_1, \lambda_2, \lambda_3, \lambda_4, \lambda_5,
\lambda_6, \lambda_7, \lambda_8, \lambda_9, \lambda_{10}$ in $[0,
0.032]$ such that
$$
\sum_{k=1}^{10}a_k g_5(\lambda_k) < 0,
$$
where
$$
a_k = \frac{p(-\lambda_k)}{\lambda_k\pi'(-\lambda_k)}
$$
and
$$
\pi(x) = \Pi_{j=1}^{10}(x + \lambda_j).
$$

We take $p(x) = x^6$ and $\lambda_k = \frac{k}{1250}$ for $1 \leq
k \leq 10$. Then by the argument as in the case of $n = 2, 3, 4$
we have
$$
\sum_{k=1}^{10}a_kg_5(\lambda_k)
= - \frac{33848952554021}{3845214843750000}  < 0,
$$
hence $g_5 \not\in M_5([0,\alpha_5])$.
\hfill\qed

\section{Comments}
Motivated by results on operator monotone and matrix
monotone functions and their relation to $C^*$-algebras
\cite{JiTomiyama}, \cite{HansenJiTomiyama-art},
 \cite{Ogasawara}, \cite{OsakaSilvestrovTomiyamaIJM-art},
\cite{Pedersenbook}, \cite{Wu-art}, \cite{SparrMathScandart}, 
and the monotonicity gap inclusion results and the $C^*$-algebraic version of interpolation spaces obtained in \cite{AKS}, we feel that the related problem of finding a $C^*$-algebraic interpretation and perhaps a $C^*$-algebraic generalization of the spaces ${\cal M}_n$ would be of interest.


\section{appendix}
In this section we put graphs which are used in the above examples.

\begin{figure}[htbp]
\begin{center}
\includegraphics[width=6cm]{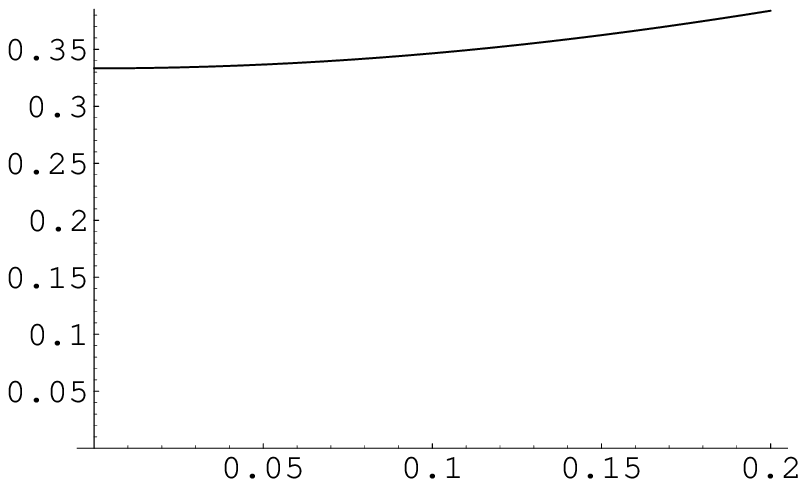}
\hspace{15mm}
\includegraphics[width=6cm]{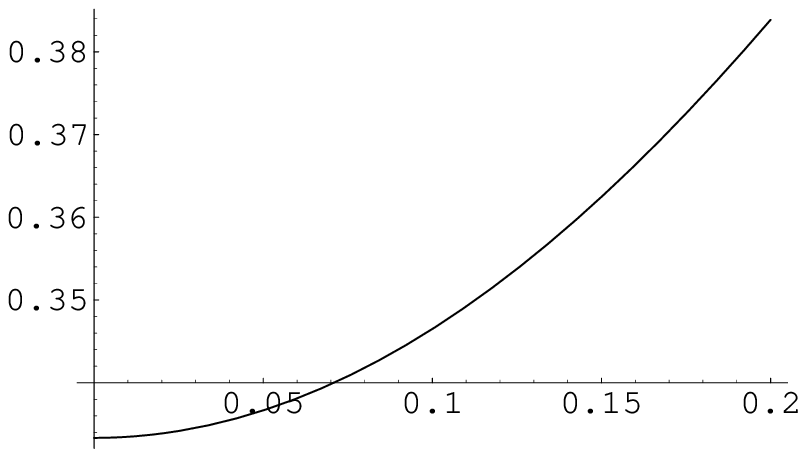}
\end{center}
\caption{(3.1) $ det(M_3(g_3;x)_{22}) $}
\end{figure}

\begin{figure}[htbp]
\begin{center}
\includegraphics[width=6cm]{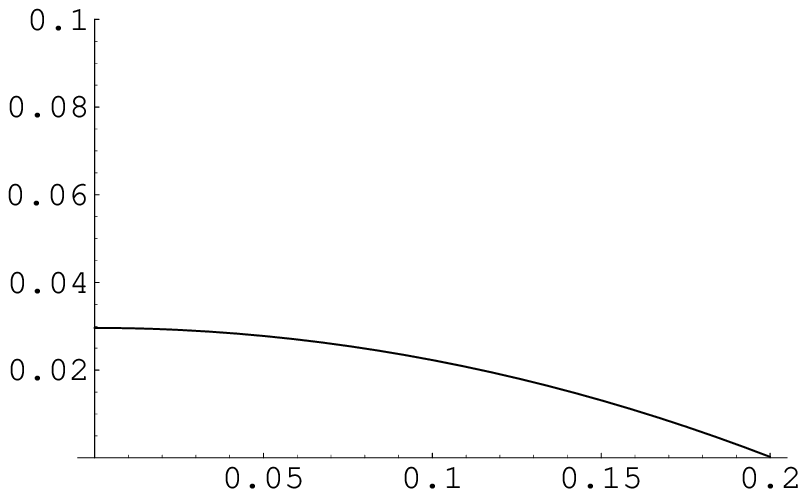}
\hspace{15mm}
\includegraphics[width=6cm]{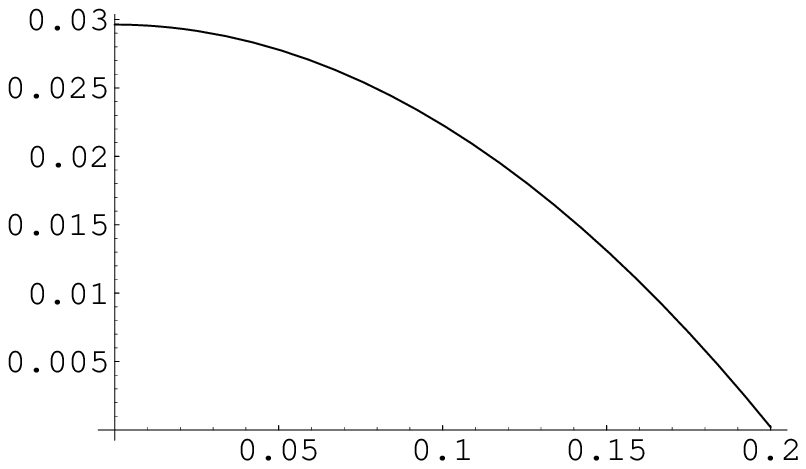}
\end{center}
\caption{(3.2) $ det(M_3(g_3;x)) $}
\end{figure}

\begin{figure}[htbp]
\begin{center}
\includegraphics[width=6cm]{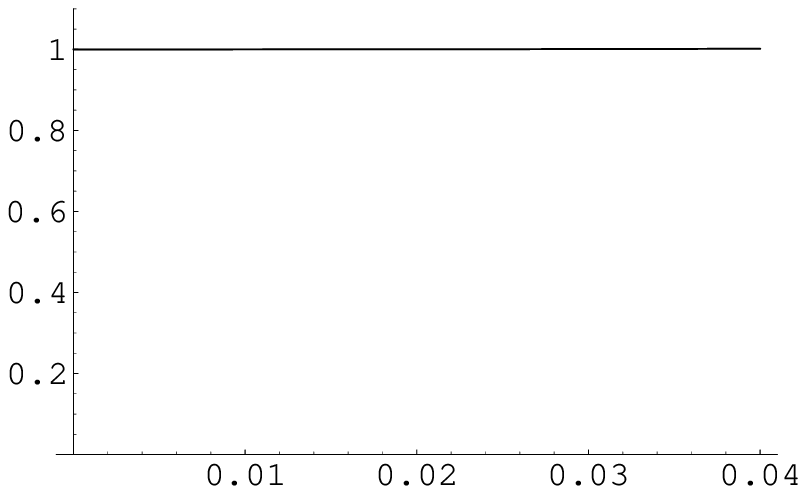}
\hspace{15mm}
\includegraphics[width=6cm]{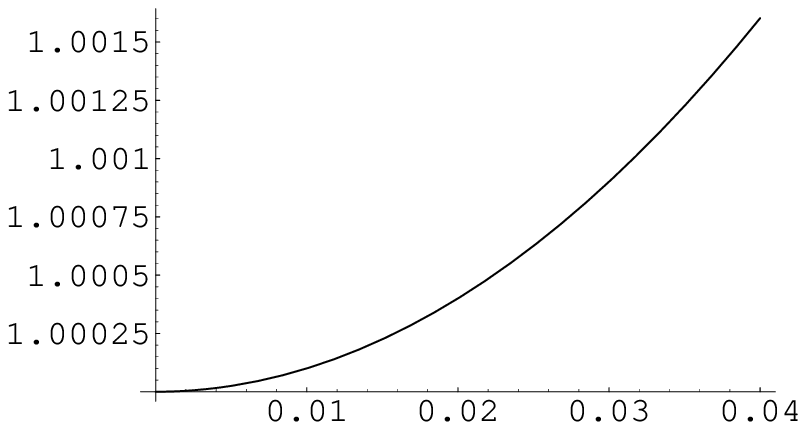}
\end{center}
\caption{(4.1) $ det(M_4(g_4;x)_{11}) $}
\end{figure}

\begin{figure}[htbp]
\begin{center}
\includegraphics[width=6cm]{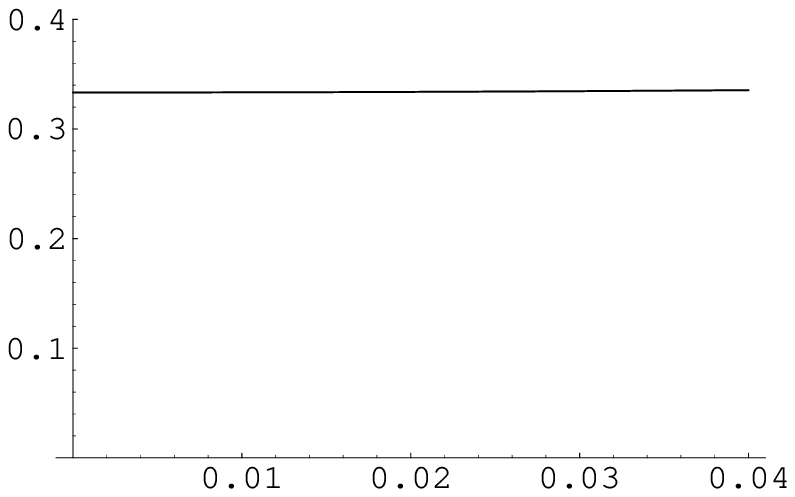}
\hspace{15mm}
\includegraphics[width=6cm]{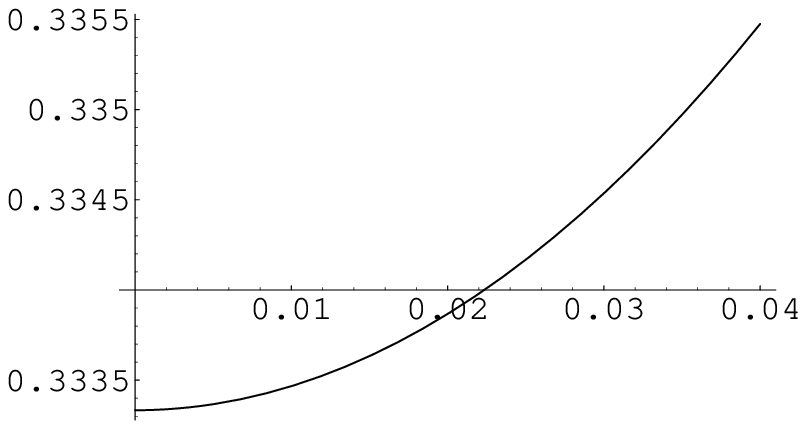}
\end{center}
\caption{(4.2) $ det(M_4(g_4;x)_{22}) $}
\end{figure}

\begin{figure}[htbp]
\begin{center}
\includegraphics[width=6cm]{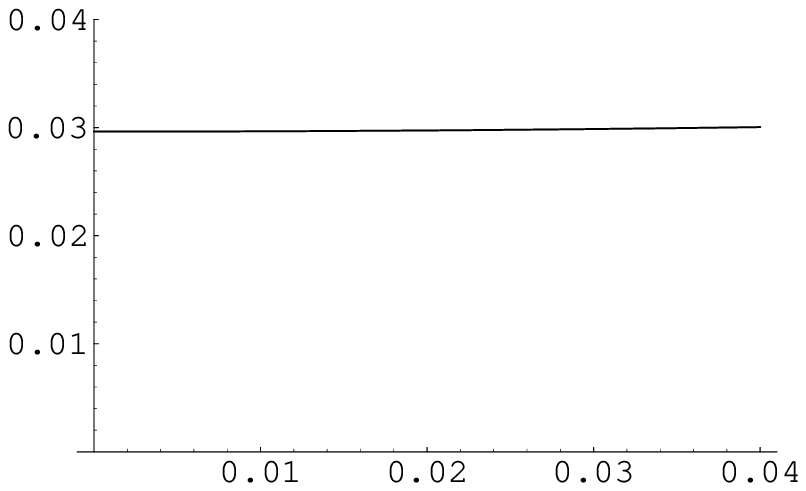}
\hspace{15mm}
\includegraphics[width=6cm]{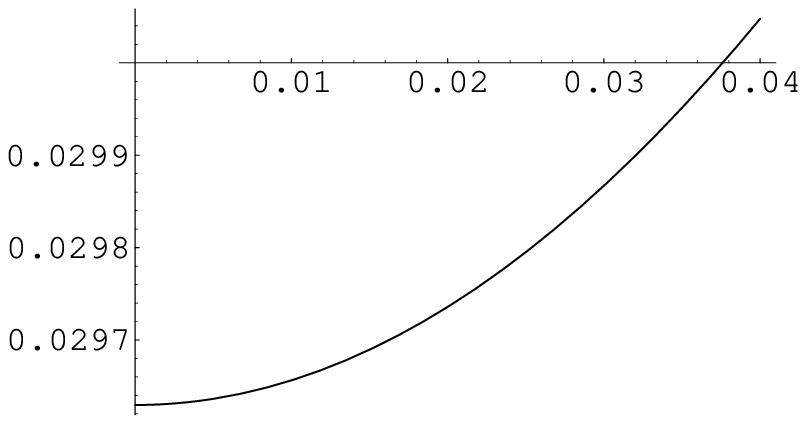}
\end{center}
\caption{(4.3) $ det(M_4(g_4;x)_{33}) $}
\end{figure}

\begin{figure}[htbp]
\begin{center}
\includegraphics[width=6cm]{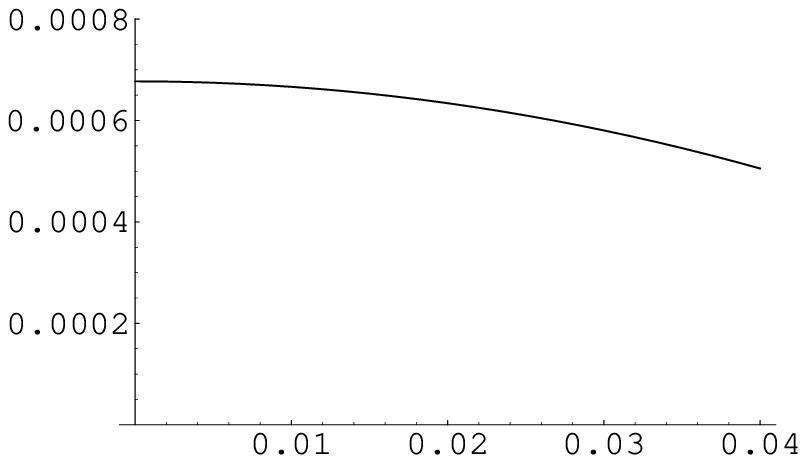}
\hspace{15mm}
\includegraphics[width=6cm]{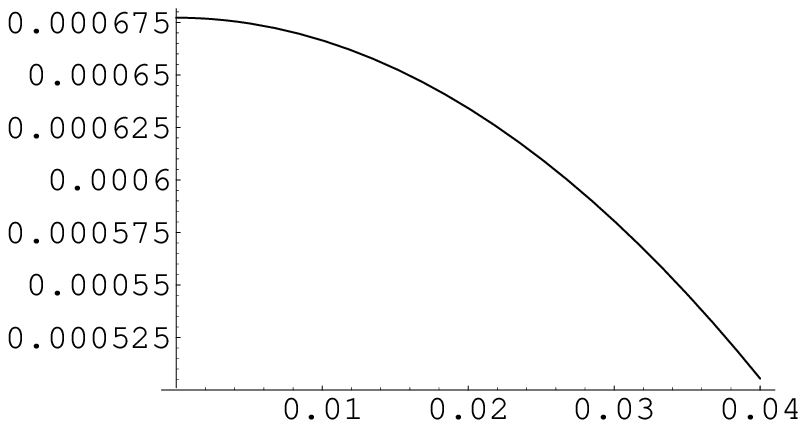}
\end{center}
\caption{(4.4) $ det(M_4(g_4;x)) $}
\end{figure}

\begin{figure}[htbp]
\begin{center}
\includegraphics[width=6cm]{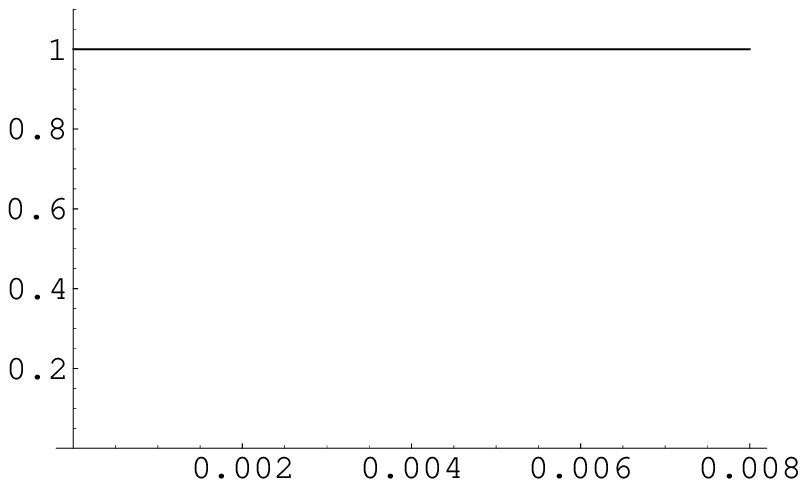}
\hspace{15mm}
\includegraphics[width=6cm]{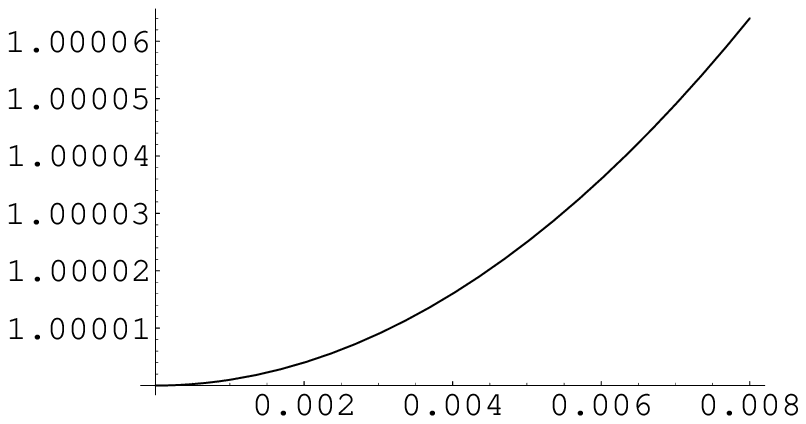}
\end{center}
\caption{(5.1) $ det(M_5(g_5;x)_{11}) $}
\end{figure}

\begin{figure}[htbp]
\begin{center}
\includegraphics[width=6cm]{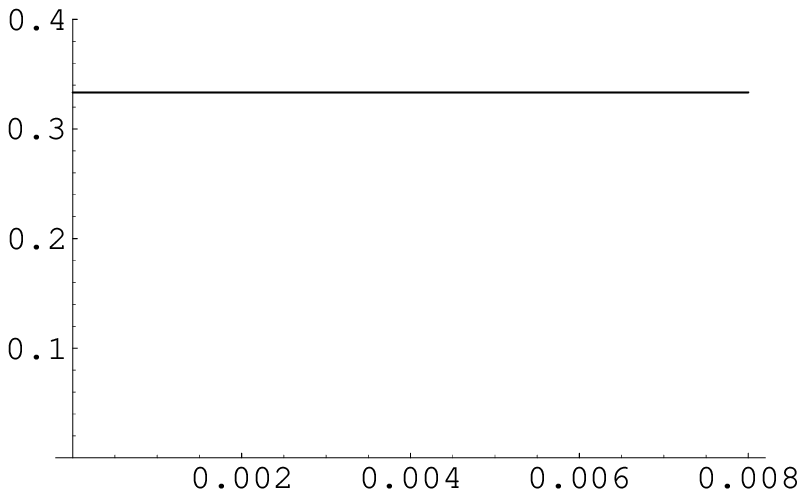}
\hspace{15mm}
\includegraphics[width=6cm]{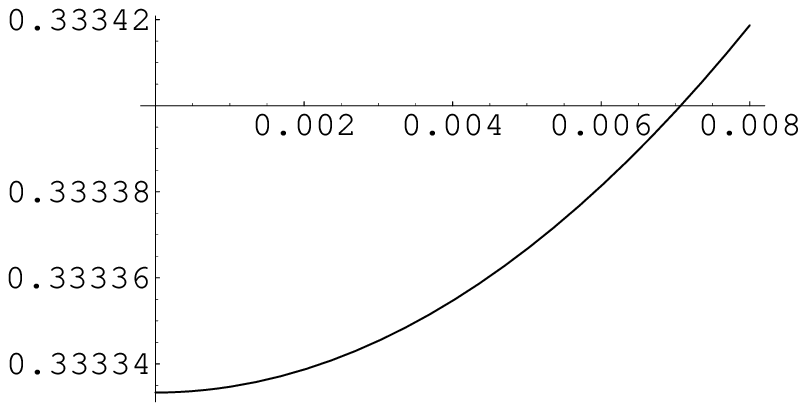}
\end{center}
\caption{(5.2) $ det(M_5(g_5;x)_{22}) $}
\end{figure}

\begin{figure}[htbp]
\begin{center}
\includegraphics[width=6cm]{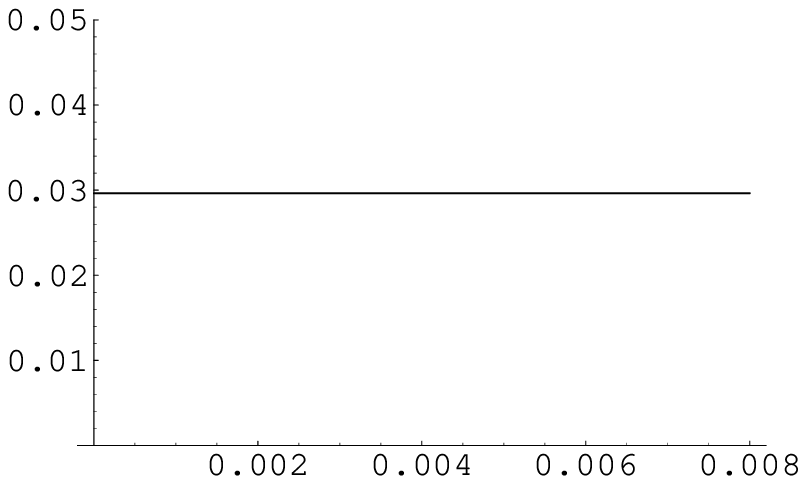}
\hspace{15mm}
\includegraphics[width=6cm]{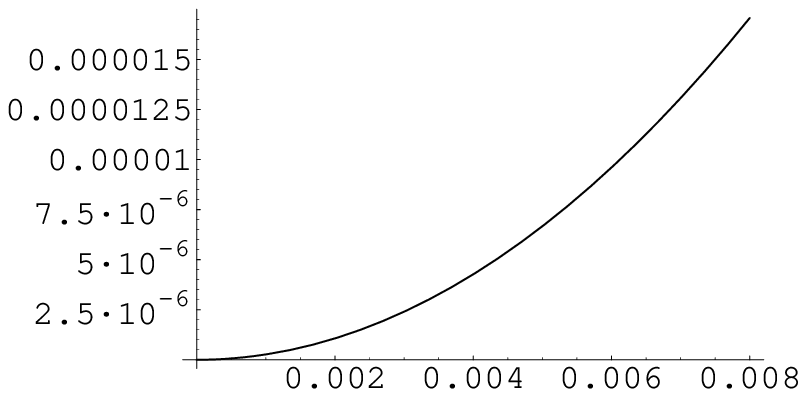}
\end{center}
\caption{(5.3) $ det(M_5(g_5;x)_{33}) $}
\end{figure}

\begin{figure}[htbp]
\begin{center}
\includegraphics[width=6cm]{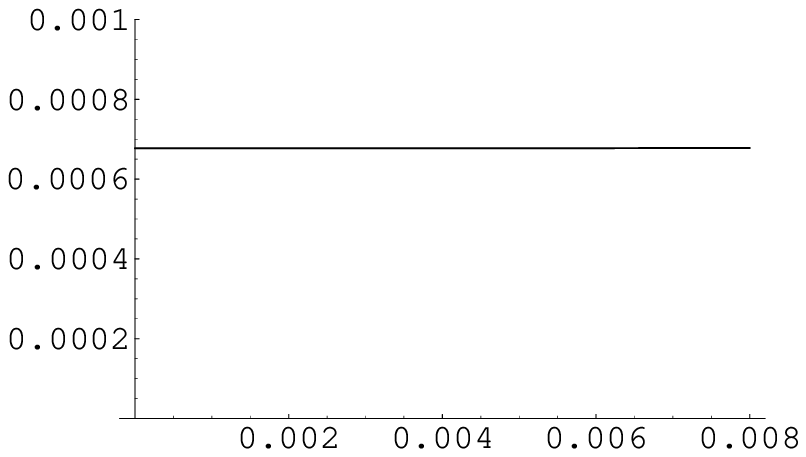}
\hspace{15mm}
\includegraphics[width=6cm]{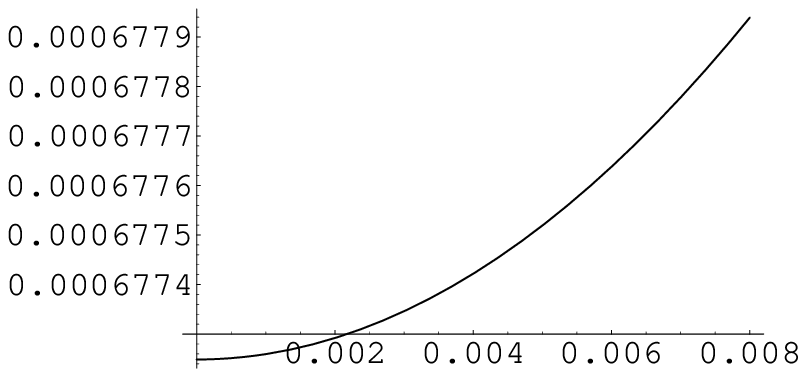}
\end{center}
\caption{(5.4) $ det(M_5(g_5;x)_{44}) $}
\end{figure}

\begin{figure}[htbp]
\begin{center}
\includegraphics[width=6cm]{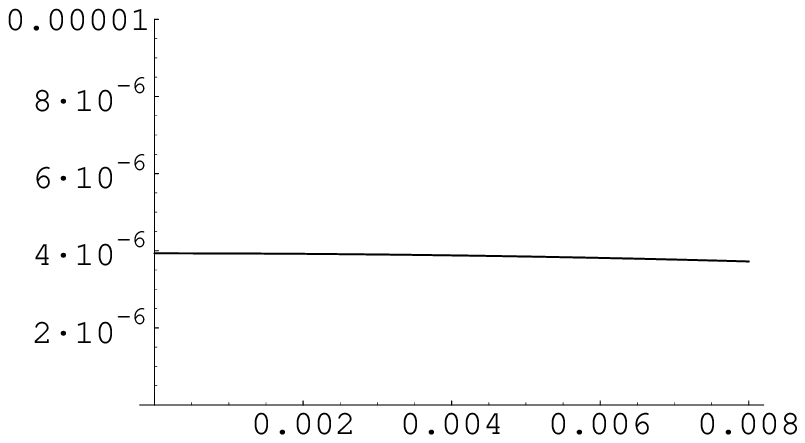}
\hspace{15mm}
\includegraphics[width=6cm]{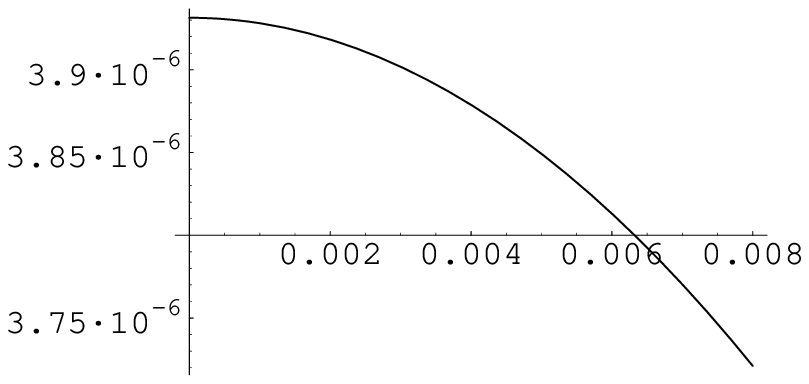}
\end{center}
\caption{(5.5) $ det(M_5(g_5;x)) $}
\end{figure}

\end{document}